\documentclass[11pt]{article}
\usepackage{amsfonts}
\usepackage{amsmath}
\usepackage{amssymb}
\usepackage{url}
\usepackage{setspace}
\usepackage{fancyheadings}

\newtheorem{theorem}{Theorem}[section]
\newtheorem{lemma}[theorem]{Lemma}
\newtheorem{proposition}[theorem]{Proposition}
\newtheorem{corollary}[theorem]{Corollary}

\newenvironment{proof}[1][Proof]{\begin{trivlist}
\item[\hskip \labelsep {\bfseries #1}]}{\end{trivlist}}
\newenvironment{definition}[1][Definition]{\begin{trivlist}
\item[\hskip \labelsep {\bfseries #1}]}{\end{trivlist}}

\newenvironment{comment}[1][Comment]{\begin{trivlist}
\item[\hskip \labelsep {\bfseries #1}]}{\end{trivlist}}

\newcommand{\ssp}{\doublespacing}

\ssp

\newcommand{\qed}{\nobreak \ifvmode \relax \else
      \ifdim\lastskip<1.5em \hskip-\lastskip
      \hskip1.5em plus0em minus0.5em \fi \nobreak
      \vrule height0.75em width0.5em depth0.25em\fi}

\begin{document}



\phantom{xxxx}
{}
\vskip 2in

\begin{center}
\Large{Cayley graphs formed by conjugate generating sets of $S_n$}\\
\large{Jacob Steinhardt}\\
\large{Thomas Jefferson High School for Science and Technology}\\
\large{\url{jacob.steinhardt@gmail.com}}\\
\large{2538 Oak Valley Drive}\\
\large{Vienna, VA, USA 22181}\\
\large{(703)-242-7397}
\end{center}

\vskip 12in




\begin{abstract}
We investigate subsets of the symmetric group with structure similar to that of a graph. The ``trees'' of these subsets correspond to minimal conjugate generating sets of the symmetric group. There are two main theorems in this paper. The first is a characterization of minimal conjugate generating sets of $S_n$. The second is a generalization of a result due to Feng characterizing the automorphism groups of the Cayley graphs formed by certain generating sets composed of cycles. We compute the full automorphism groups subject to a weak condition and conjecture that the characterization still holds without the condition. We also present some computational results in relation to Hamiltonicity of Cayley graphs, including a generalization of the work on quasi-hamiltonicity by Gutin and Yeo to undirected graphs.
\end{abstract}


\vskip 12in


\section{Terminology}

In this paper, we will let $N$ denote the set $\{1,2,\ldots,n\}$. $S_n$ will denote the symmetric group acting on $n$ elements with canonical action on $N$. $A_n$ will denote the alternating group acting on $N$. We will use the notation $(a_1 \ a_2 \ \ldots \ a_{k_1})(b_1 \ b_2 \ \ldots \ b_{k_2})\ldots$ to express a permutation as a product of disjoint cycles. By the \textit{support} of a permutation we will mean those elements not fixed by the permutation. Given permutations $\sigma$, $\tau$, $\sigma\tau$ denotes $\tau \circ \sigma$.

Given a multiset $A = (a_1,\ldots,a_k)$ with $a_i \geq 2$, we define its \textit{extended conjugacy class in $S_n$} to be the set of all permutations such that, when decomposed into disjoint cycles, contain cycles of lengths $a_1,\ldots,a_k$, and no others. We denote it by $\mathcal{C}(A)$.

Given a set $S \subset N$, define the \textit{subsymmetric group of $S$} as the set of all permutations in $S_n$ that fix all elements outside of $S$. Define the \textit{subalternating group of $S$} as the set of all even permutations in $S_n$ (i.e., the permutations in $A_n$) that fix all elements outside of $S$. A \textit{semisymmetric group of $S$} is defined as a subgroup of $S_n$ that stabilizes $S$ whose restriction to $S$ forms a symmetric group acting on $S$. A \textit{semialternating group of $S$} is defined similarly.

Given a graph $\Gamma$, we let $V(\Gamma)$ and $E(\Gamma)$ denote the vertices and edges of $\Gamma$, respectively. An \textit{Eulerian path} is a walk in $\Gamma$ that traverses each edge exactly once. It is called an \textit{Eulerian cycle} if the first and last vertices in the walk are the same. Given a group $G$ and a set $S \subset G$, the \textit{Cayley graph} $\Gamma = Cay(G,S)$ is defined as follows: each vertex is an element of $G$, and two vertices $g,h \in V(\Gamma)$ are adjacent if $gh^{-1} \in S$ or $hg^{-1} \in S$.

\section{Motivation and Overview}

Cayley graphs are of general interest in the field of Algebraic Graph Theory and also have certain properties desirable in practical applications. We present here a brief survey of some of the broader results and conjectures surrounding Cayley graphs. Godsil and Royle \cite{Roy} provide a useful overview of work on graphs with transitive permutation groups in general, which we partially reproduce here. First, all Cayley graphs are vertex-transitive since the mapping $\phi_g(x) = xg$ is an automorphism for all $g \in G$. As such, there is always a representation of $G$ in $Aut(Cay(G,S))$, denoted $R(G)$. $R(G)$ acts not only transitively but regularly on the vertices of $Cay(G,S)$. Sabidussi has shown that the converse of this is true, namely that $\Gamma$ if a Cayley graph of $G$ if and only if $Aut(\Gamma)$ contains a subgroup isomorphic to $G$ that acts regularly on $V(\Gamma)$ \cite{Sab}.

Minimally generated Cayley graphs have provably maximal vertex connectivity \cite{Roy}, which points to uses in practical applications. Specifically, Cayley graphs have been used to create networks with small diameter and valency and high connectivity for uses in parallel processing, and Schreier coset graphs, a generalization of Cayley graphs, have been used to solve certain routing problems \cite{Coo1}. See also \cite{Coo2} and \cite{Lew}.

One tempting conjecture related to Cayley graphs is that every Cayley graph has a Hamiltonian cycle, first given by Strasser in $1959$ \cite{Rap}. For more information, see \cite{Dob}, \cite{Kut}, and \cite{Pak}. While we offer some computational ideas in relation to this conjecture, the focus of this paper is on another difficult problem in Algebraic Graph Theory, that of characterizing the automorphism groups of Cayley graphs. We have a poor understanding of the automorphism groups of Cayley graphs, though there are some notable exceptions (see below). These groups are fundamental as the most natural algebraic structure to associate with an arbitrary highly symmetric graph.

Note that a graph can be defined as a collection of vertices and edges. Two vertices are adjacent if there exists an edge connecting them, and two vertices $v_1$ and $v_2$ are connected if there exists a sequence of adjacent vertices containing $v_1$ and $v_2$. On the other hand, consider the following definition: Given a collection of vertices $V$ and a collection of edges $E$, we can let each element of $E$ act on $V$ as a transposition swapping the two vertices on which $E$ is incident. If we then let multiplication in $E$ extend through the definitions of a group action, $E$ generates a subgroup of the symmetric group acting on $V$ (we denote this subgroup as $<E>$). Then we say that $v_1$ and $v_2$ are adjacent if $(v_1 \ v_2) \in E$, and that $v_1$ and $v_2$ are connected if $(v_1 \ v_2) \in <E>$. Additionally, connected components correspond to orbits of $V$ under $E$. A tree is a minimal generating set of $S_{|V|}$ consisting only of transpositions (thus the fact that trees have $n-1$ vertices corresponds to the fact that it takes $n-1$ transpositions to generate $S_n$). It is easily verified that these definitions are equivalent.

The algebraic properties of the related Cayley graphs of trees in the above definition are well-understood. We know in particular that these graphs are Hamiltonian \cite{Kom}. Furthermore, in $2003$ Feng \cite{Fen2} generalized a result by Godsil \cite{Roy} that fully characterizes the automorphism groups of these graphs.

\section{$\mathcal{C}$-trees}

\subsection{Definitions}

The above definition of a graph in terms of transpositions can be generalized. Given a collection of vertices (which, from now on, for convenience, will without loss of generality be $N$), and a set $T \subset S_n$ in which all elements of $T$ are conjugate (say with conjugacy class $\mathcal{C}$), then we can define elementary notions in a $\mathcal{C}$-graph as follows. $v_1,v_2 \in N$ are adjacent if they have the same orbit under a single element of $T$. They are \textit{semi-connected} if they have the same orbit under $T$, and \textit{connected} if $(v_1 \ v_2) \in <T>$ (it is then easy to verify that semi-connectivity and connectivity are equivalence relations). Connected components correspond to subsymmetric groups of $<T>$. A \textit{tree} is a minimal generating set of $S_n$ with all elements lying in $\mathcal{C}$. It is natural to ask why we add the somewhat artificial-looking stipulation that all elements of $T$ belong to the same conjugacy class. The main reason is that this stipulation is inherent in the construction of a normal graph, where all edges are transpositions. Additionally, without this restriction we get the result that a tree, under our fairly intuitive definition, almost always has $2$ edges since $(1 \ 2)$ and $(2 \ 3 \ 4 \ \ldots \ n)$ generate $S_n$.

In this paper we will characterize $\mathcal{C}$-trees and study some of their properties, including a generalization of Feng's result. However, we will still use the language of graphs for the sake of intuition. For approaches to extending the above intuitive generalization to a well-structured system, see the concluding section on open problems.

With $\mathcal{C}$-trees defined, we introduce some more terminology associated to them. A set $T \subset S_n$ is said to be \textit{semi-connected} if $N$ has a single orbit under $T$ (i.e. all elements of $N$ are semi-connected). We call it \textit{split} if the intersection of the supports of any two elements of $T$ has size at most one. Note that if $T$ generates $S_n$ then it must be connected.

Given a multiset $A$ and integer $n$, we define $f(A,n)$ to be the infimum of $|G|$ across all $G \subset \mathcal{C}(A)$ that generate $S_n$. We aim to find $f(A,n)$ for every $A$ for sufficiently large $n$. Let $c(A)$ be defined as

$$\sum_{i=1}^{|A|} a_i-1$$

We aim to prove that there exists a function $X_0(A)$ such that, for $n \geq X_0(A)$, $f(A,n)$ is equal to

\begin{equation}
\label{SizeBound}
\lceil\frac{n-1}{c(A)}\rceil
\end{equation}

\noindent when $c(A)$ is odd, and $\infty$ otherwise. When $c(A)$ is odd, then $A$ defines the conjugacy class of an even permutation and so $f(A,n)$ is obviously $\infty$ (because it is impossible to generate any odd permutations). Note further that $f(A,n)$ is necessarily at least the value given by (\ref{SizeBound}), as $c(A)$ counts the number of transpositions necessary to generate an element of $\mathcal{C}(A)$, and so if it was smaller then it would be possible to generate $S_n$ with less than $n-1$ transpositions. Another way to see this is that no potential generating sets can be semi-connected, and thus cannot generate $S_n$.

\subsection{Cycles}

We study first the case of a single $k$-cycle, i.e. $|A| = 1$ and $a_1=k$. We will give explicit generators for $S_{2k-1}$:

\begin{proposition}
\label{CycleTreeSize}
The set $\{(1 \ 2 \ \ldots \ k),(k \ k+1 \ \ldots \ 2k-1)\}$ generates $S_{2k-1}$ when $k$ is even.
\end{proposition}

\begin{proof}
We construct something similar to a semisymmetric group of $\{1,2,\ldots,k\}$, except with the elements lying in the positions $\{k,k+1,\ldots,2k-1\}$. From this we will construct a subsymmetric group of $\{1,2,\ldots,k\}$, which will finally allow us to construct the entire symmetric group $S_{2k-1}$.

\begin{lemma}
\label{SemiSymmetricCycle}
We can place the elements $\{1,2,\ldots,k\}$ in any order in positions $\{k,k+1,\ldots,2k-1\}$ as long as we allow the other elements to move arbitrarily (even for odd $k$).
\end{lemma}

\begin{proof}
Let $\sigma = (1,2,\ldots,k)$ and $\tau = (k,k+1,\ldots,2k-1)$. If we allow the first $k-1$ elements to be arbitrary, we can set the last $k$ elements to any permutation $\pi = (\pi_1,\pi_2,\ldots,\pi_{k})$ of $(1,2,\ldots,k)$ as follows: Rotate $\pi_k$ to the $k$th position using $\sigma$, then apply $\tau$ once. Now rotate $\pi_{k-1}$ to the $k$th position (again with $\sigma$), and apply $\tau$ again. Continue this process until we have put all $k$ of the desired elements into place. For example, to make the last $4$ elements $(2,4,1,3\}$, we would apply $\sigma^1\tau\sigma^2\tau\sigma^1\tau\sigma^2$. Each successive set of applications of $\sigma$ moves the next desired element into place ($3$, then $1$, then $4$, then $2$). $\qed$
\end{proof}

\begin{lemma}
\label{CyclesGenerate}
For $k$ even, the $k$-cycles generate $S_n$ for $n \geq k$. Moreover, for $k$ odd, the $k$-cycles generate $A_n$.
\end{lemma}

\begin{proof}
Let $U$ be the set of all $k$-cycles. $U$ is closed under conjugation, so $<U>$ is normal in $S_n$. Thus $<U>$ is either $(e)$, $D_4$, $A_n$, or $S_n$, where by $D_4$ we mean the normal subgroup of $A_4$ isomorphic to the dihedral group of order $4$. It can't be $(e)$ because $U$ is non-trivial, and it can't be $D_4$ because $D_4$ contains no cycles. Thus it is $A_n$ or $S_n$. If $k$ is odd, it must be $A_n$ because $U$ consists of only even permutations. If $k$ is even, it must be $S_n$ since $U$ contains an odd permutation. $\qed$

\end{proof}

\begin{lemma}
\label{SubSymmetricCycle}
We can generate the subsymmetric group on $\{1,2,\ldots,k\}$ when $k$ is even.
\end{lemma}

\begin{proof}
Take some such permutation $\pi$ generated in the manner of Lemma~\ref{SemiSymmetricCycle}, and consider $\pi\tau\pi^{-1}$. This creates an arbitrary $k$-cycle among the first $k$ elements while fixing the last $k-1$ elements. Then by Lemma~\ref{CyclesGenerate} we can generate the subsymmetric group on $\{1,2,\ldots,k\}$. $\qed$
\end{proof}

Now, to generate an arbitrary permutation $\pi = \{\pi_1,\ldots,\pi_{2k-1}\}$ in $S_{2k-1}$, first use $\tau$ to move $\{\pi_{k+1},\ldots,\pi_{2k-1}\}$ to the first $k-1$ elements of the set. We can do this by moving each one to the $k$th position, then, since we can generate any permutation among the first $k$ elements (recall that we can do this by Lemma~\ref{SubSymmetricCycle}), move it to an arbitrary place among the first $k-1$ elements in which we haven't already put anything with this process. Next apply a permutation that puts $\{\pi_{k+1},\ldots,\pi_{2k-1}\}$ in the proper order (though leaving them in the positions $\{1,2,\ldots,k-1\}$). We can then move them to their proper locations with $(\sigma\tau)^{k-1}$. Now that the last $k-1$ elements are in place, we can apply whatever permutation is necessary to put the first $k$ elements in place. We can thus generate an arbitrary permutation and therefore $S_{2k-1}$. This completes Proposition~\ref{CycleTreeSize}. $\qed$
\end{proof}

\begin{proposition}
\label{GeneralCycleTreeSize}
$f((k),n(k-1)+1) = n$ for $n \geq 2$ and $k$ even.
\end{proposition}

\begin{proof}
This follows by induction on $n$. Proposition~\ref{CycleTreeSize} proves the base case of $n=2$. The inductive step is completed by the following easily verified lemma:

\begin{lemma}
\label{SubSymmetricExtension}
The subsymmetric group on $S$, together with the cycle $(a_1,\ldots,a_n)$, generates the subsymmetric group on $S \cup \{a_1,\ldots,a_n\}$ provided that $S \cap \{a_1,\ldots,a_n\} \neq \emptyset$ and $S \not\subset \{a_1,\ldots,a_n\}$.
\end{lemma}
\end{proof}

\begin{corollary}
\label{FinalCycleTreeSize}
$f((k),n) = \lceil\frac{n-1}{k-1}\rceil$ for $n \geq 2k-1$.
\end{corollary}

\begin{proof}
Take the construction for when $\frac{n-1}{k-1}$ is an integer (i.e. that given above in Proposition~\ref{GeneralCycleTreeSize}). Then, to extend the formula to arbitrary $n$, add the $k$-cycle $(n-k+1,n-k+2,\ldots,n)$ and apply Lemma~\ref{SubSymmetricExtension}. This completes the claimed characterization of $f(A,n)$ for cycles.
\end{proof}

\subsection{Products of Transpositions}

Having proven our result for cycles, we would like to extend it to more complex permutations. We will start with the simplest of these, i.e. products of disjoint transpositions. We call a permutation \emph{basic} if it is of this form, and denote $B_k = (2,2,\ldots,2)$ ($k$ twos).

\begin{proposition}
\label{BasicTreeSize}
$f(B_k,n) = \lceil\frac{n-1}{k}\rceil$ for $n \geq k(2k+1)+1$ and $k$ odd.
\end{proposition}

\begin{proof}
Our goal will be to write $2k+1$ explicit generators for $S_{k(2k+1)+1}$. Then we can easily proceed by induction as before. Such generators must form a semi-connected set. However, we would also like the set to be split so that only one cycle interacts at a time when multiplying permutations. To do this, we will find an Eulerian cycle of $K_{2k+1}$, which will allow us to create a connected and split set.

First note that every vertex of $K_{2n+1}$ has even degree (in particular, degree $2n$), so that $K_{2n+1}$ has an Eulerian cycle. We now construct generators from this cycle. We start with an example, then give a general method. Consider $k=3$, $2k+1=7$, and the cycle $1 \to 2 \to 4 \to 7 \to 1 \to 3 \to 6 \to 7 \to 2 \to 5 \to 6 \to 1 \to 4 \to 5 \to 7 \to 3 \to 4 \to 6 \to 2 \to 3 \to 5 \to 1$. We have the generators

\singlespacing
\[ g_1 = (1 \ 2)(5 \ 6)(12 \ 13) \ \ g_2 = (2 \ 3)(9 \ 10)(19 \ 20) \]
\[ g_3 = (6 \ 7)(16 \ 17)(20 \ 21) \ \ g_4 = (3 \ 4)(13 \ 14)(17 \ 18) \]
\[ g_5 = (10 \ 11)(14 \ 15)(21 \ 22) \ \ g_6 = (7 \ 8)(11 \ 12)(18 \ 19) \ \ g_7 = (4 \ 5)(8 \ 9)(15 \ 16) \]
\\
\noindent \ssp Note that if we follow the path of permutations containing $(1,2),(2,3),(3,4),\ldots,(21,22)$, then we get $g_1,g_2,g_4,g_7,\ldots,g_3,g_5$, i.e. exactly the constructed Eulerian cycle (with the exception of the final vertex). This is the general method in which we will construct our generators. Specifically, we place the transposition $(i \ i+1)$ in the generator corresponding to the $i$th vertex visited in the cycle. Note that the properties of an Eulerian cycle guarantee that these generators will be semi-connected and split. The semi-connected part is obvious, whereas the split part is a consequence of the fact that every edge is traversed exactly once, which corresponds to the fact that each pair of generators move at most one common element. We now show this construction generates $S_{k(2k+1)+1}$.

\begin{theorem}
\label{ConnectedSplitTranspositions}
If $T \subset \mathcal{C}(B_k)$ is semi-connected and split, then $<T> = S_n$ or $A_n$, depending on whether $k$ is even or odd.
\end{theorem}

\begin{proof}
We call two generators \textit{adjacent} if they move a common element. Consider two adjacent generators, $g_i$ and $g_j$, and consider $g_ig_jg_ig_j$. All transpositions are applied twice in this case and therefore cancel, except for the two transpositions that act on the same element, which multiply to a three cycle. So, in the above example, $g_5g_7g_5g_7 = (14 \ 15)(15 \ 16)(14 \ 15)(15 \ 16) = (14 \ 16 \ 15)^2 = (14 \ 15 \ 16)$. In this manner, we generate all $3$-cycles of the form $(i \ i+1 \ i+2)$. We wish to show that these generate $A_n$. From this, we would be done, since any odd permutation then allows us to generate $S_n$. In fact, it is convenient for later purposes to prove a slightly stronger result:

\begin{lemma}
\label{SubAlternatingExtension}
When $n$ is odd, the subalternating group on $S$, together with the cycle $(a_1,\ldots,a_n)$, generates the subalternating group on $S \cup \{a_1,\ldots,a_n\}$ provided that $S \cap \{a_1,\ldots,a_n\} \neq \emptyset$ and that $|S \cap \{a_1,\ldots,a_n\}| \leq |S|-2$. When $n$ is even, it generates the entire subsymmetric group.
\end{lemma}

\begin{proof}
Like Lemma~\ref{SubSymmetricExtension}, the proof is easy enough to omit. The only important detail to note is that we get $\frac{|S \cap \{a_1,\ldots,a_n\}|!}{2}$ distinct permutations, which must be the alternating group when $n$ is odd or must generate the symmetric group by Lagrange's theorem when $n$ is even.
\end{proof}

In particular, a $3$-cycle looks like $A_3$, so the given $3$-cycles indeed generate the alternating group (they are semi-connected since $T$ was semi-connected), and we are done with Theorem~\ref{ConnectedSplitTranspositions}.
\end{proof}

Our proof of the remainder of Proposition~\ref{BasicTreeSize} (i.e. the induction and extension to cases when $c(A)$ does not divide the $n-1$) follows in exactly the same manner as that of Proposition~\ref{CycleTreeSize}, and so we omit it, instead referring the reader to Proposition~\ref{GeneralCycleTreeSize} and Corollary~\ref{FinalCycleTreeSize}. The only necessary modification is that we must deal with each of the cycles in the added permutation one at a time in our inductive step.
\end{proof}

\subsection{The General Case}

We would next like a general criterion for connectedness. We present it here:

\begin{definition}
A set $T \subset \mathcal{C}(A)$ is called \textit{balanced} if it is possible to divide the set of orbits of elements of $T$ into disjoint sets $S_1,S_2,\ldots,S_{|A|}$ such that all orbits in $S_i$ have the same size and each element of $S_i$ overlaps with at least one other element of $S_i$. We will denote the size of the orbits in $S_i$ by $a_i$.
\end{definition}

\begin{theorem}
\label{ConnectedSplitBalanced}
All semi-connected, split, balanced sets in the extended conjugacy class of an odd permutation generate $S_n$.
\end{theorem}

\begin{proof}
We proceed by induction on two quantities: first $|A|$, then the number of occurrences of $2$ in $A$. Note that we have already proven the base cases of this induction in Corollary~\ref{FinalCycleTreeSize} and Proposition~\ref{BasicTreeSize}.

We call two permutations \textit{$i$-adjacent} if they both have orbits in $S_i$. Pick $i$ such that $a_i$ is maximal, and consider any two $i$-adjacent permutations, $\sigma$ and $\tau$, with orbits $_i\sigma$ and $_i\tau$ in $S_i$. By the same argument as Proposition~\ref{CycleTreeSize}, these generate the semialternating group on the elements moved by the two identified cycles (moving $2a_i-1$ elements in total). Thus in particular, by \textbf{Chebyshev's Theorem} \cite{Ram}, there exists a prime strictly between $a_i$ and $2a_i$, and the semialternating group contains a cycle of this length, call it $p$.

Consider each cycle of this length in our semialternating group, and apply it $\displaystyle lcm_{a \in A} a$ times. Since $p$ is prime and $a_j < p$ for each $j$, we end up with a $p$-cycle, which we can then apply some number of times to get back to our original $p$-cycle. Note, however, that all other elements that were moved contained only cycles of length $a_j$ for some $j$, and so were all cancelled out by the above repeated application. Thus we are left only with the actual $p$-cycle. Repeating this for all such $p$-cycles in the semialternating group gives us all actual $p$-cycles, i.e. those living in the associated subalternating group. Thus, by the same arguments as in Lemma~\ref{CyclesGenerate}, they generate the entire subalternating group. If $a_i$ is odd, then $_i\sigma$ and $_i\tau$ both live in this subalternating group, and so we can take $\sigma(_i\sigma)^{-1}$ and $\tau(_i\tau)^{-1}$. Taking $\sigma(_i\sigma)^{-1}$ for all $\sigma \in T$ (we can do this since $T$ is balanced) gives us a semi-connected, split, balanced set with strictly less orbits in each permutation, so that we can apply the inductive step to generate the subsymmetric group on the elements moved by these new permutations. Then, by adding $_i\sigma$ for each $\sigma \in T$, by Lemma~\ref{GeneralCycleTreeSize} we can generate the entire symmetric group, and we are done.

On the other hand, if $a_i$ is even, then we can only cancel $_i\sigma$ down to a transposition. However, this gives us an extended conjugacy class with the same number of orbits, but with strictly more occurrences of $2$ in $A$ than before. Thus we can apply the inductive step in the same manner as above, and are once again done. Note that in both cases we attained elements in the new extended conjugacy class by multiplying elements in the old extended conjugacy class by an even permutation. This shows that the new extended conjugacy class does indeed correspond to an odd permutation. $\qed$
\end{proof}

We are now ready to prove our major contention:

\begin{theorem}
\label{GeneralTreeSize}
Let $A$ be a multiset of size $k$. Then there exists some $X_0(A)$ such that $f(A,n) = \lceil\frac{n-1}{c(A)}\rceil$ when $n \geq X_0$. Furthermore, $X_0((k)) \leq 2k-1$, $X_0(B_k) \leq \binom{2k+1}{2}+1$, and $X_0(A) \leq c(A)\Phi_{2|A|}(2|A|)+1$, where $\Phi_k$ denotes the $k$th cyclotomic polynomial.
\end{theorem}

\begin{proof}
Note that the first two bounds have already been proven. For the final case, we again use Eulerian cycles, this time with the goal of creating a semi-connected, split, and balanced set. In particular, we find a prime congruent to $1$ mod $2k$. We know that such a prime exists that is less than $\Phi_{2k}(2k)$ (proof sketch in appendix).

Take such a prime, $p$, fitting the properties described above. If $\frac{p-1}{2} = kn$, then we will work in the extended conjugacy class that is equivalent to $n$ copies of $A$, then use this to move down to $A$ itself. We construct an Eulerian cycle for $K_p$ as follows. The edges (mod $p$) will be 

$$1,2,\ldots,p,2,4,\ldots,2p,3,\ldots,3p,\ldots,\frac{p-1}{2},p-1,\ldots,\frac{p(p-1)}{2}$$

So, for example, if $p=7$ then we have (in the case of $2$-cycles) the associated generators

\singlespacing
\[ (1 \ 2)(11 \ 12)(19 \ 20) \ \ (2 \ 3)(8 \ 9)(17 \ 18) \]
\[ (3 \ 4)(12 \ 13)(15 \ 16) \ \ (4 \ 5)(9 \ 10)(20 \ 21) \]
\[ (5 \ 6)(13 \ 14)(18 \ 19) \ \ (6 \ 7)(10 \ 11)(16 \ 17) \ \ (7 \ 8)(14 \ 15)(21 \ 22) \]
\\
\noindent \ssp We can extend this past $2$-cycles (for example, permutations in the extended conjugacy class $(2,4,5)$ in the following manner:

\singlespacing
\[ (1 \ 2)(11 \ 23 \ 24 \ 12)(19 \ 37 \ 38 \ 39 \ 20) \ \ (2 \ 3)(8 \ 25 \ 26 \ 9)(17 \ 40 \ 41 \ 42 \ 18) \]
\[ (3 \ 4)(12 \ 27 \ 28 \ 13)(15 \ 43 \ 44 \ 45 \ 16) \ \ (4 \ 5)(9 \ 29 \ 30 \ 10)(20 \ 46 \ 47 \ 48 \ 21) \]
\[ (5 \ 6)(13 \ 31 \ 32 \ 14)(18 \ 49 \ 50 \ 51 \ 19) \ \ (6 \ 7)(10 \ 33 \ 34 \ 11)(16 \ 52 \ 53 \ 54 \ 17) \ \ (7 \ 8)(14 \ 35 \ 36 \ 15)(21 \ 55 \ 56 \ 57 \ 22) \]
\\
\noindent \ssp Note that we simply add elements to cycles in the $i$th column until the cycles in that column have length $a_i$. Note also that this is a balanced set by construction. It is easy to verify that this also defines an Eulerian cycle, and is thus connected and split. On the other hand, we have the following result:

\begin{lemma}
\label{Divisors}
If $B$ is equivalent to $k$ copies of $A$, and if there exists $T \subset \mathcal{C}(B)$ that generates $S_n$, then there exists $T' \subset \mathcal{C}(A)$ that generates $S_n$, and furthermore such that $|T'| = k|T|$.
\end{lemma}

\begin{proof}
Split each $\sigma \in T$ into $k$ permutations such that each of these permutations lies in $\mathcal{C}(A)$. These obviously generate $S_n$ since products of them generate $S_n$. $\qed$
\end{proof}

This proves the base case of a final induction showing that $f(A,n) = \lceil\frac{n-1}{c(A)}\rceil$ for all $n \geq X_0$, where $X_0=pc(A)+1$. This induction will finally prove Theorem~\ref{GeneralTreeSize}. However, once again this new induction is identical to the completions of Propositions~\ref{CycleTreeSize} and ~\ref{BasicTreeSize}, and so we refer the readers there for the completion of the proof.
\end{proof}

\section{Automorphism Groups}

We devote this section to the characterization of the automorphism groups of certain $\mathcal{C}$-graphs.

\begin{definition}
Given a split set of cycles $T \subset S_n$, the \textit{cycle graph} $Cyc(T)$ is formed by associating each vertex with an element of $N$ and drawing edges $x_1x_2,x_2x_3,\ldots,x_{k-1}x_k$ if $(x_1 \ x_2 \ \ldots \ x_k)$ is in $T$. Note that this involves arbitrarily choosing a ``starting'' and ``ending'' point for each cycle in $T$. When $T$ consists of transpositions, Feng \cite{Fen2} refers to $Cyc(T)$ as $Tra(T)$.
\end{definition}

\begin{definition}
Given a split set of cycles $T$, the \textit{degree} of some $t \in T$ is defined as the number of distinct points in its support that overlap with other cycles. If $t$ has degree $1$, we call it a \textit{leaf}.
\end{definition}

\begin{definition}
A split set of more than two cycles generating $S_n$ is said to be \textit{normal} if any element is adjacent to at most $1$ leaf, and furthermore $Cyc(T)$ is a tree (note that this is stronger than asking that $T$ be a minimal generating set of $S_n$, as it effectively adds the criterion that $n \equiv 1 \pmod k$, where $T$ consists of $k$-cycles).
\end{definition}

We use this to offer a partial generalization to a theorem by Feng \cite{Fen2} that states that $Aut(Cay(S_n,T)) \cong R(S_n) \rtimes Aut(S_n,T)$, where $Aut(S_n,T) = \{\phi \in Aut(S_n) \ | \ \phi(T) = T\}$, and furthermore that $Aut(S_n,T) \cong Aut(Tra(T))$. In the following, $T$ will always be normal, and if we talk about a graph it will be $Cay(S_n,T)$ unless otherwise specified:

\begin{theorem}
\label{CayleyAutomorphism}
Let $T$ be a normal set. Then $Aut(Cay(S_n,T)) \cong R(S_n) \rtimes Aut(S_n,T)$, where $R(S_n)$ is the representation of $S_n$ as an action on $Cay(S_n,T)$.
\end{theorem}

\begin{proof}
We use Feng's idea of finding cycles that force graph automorphisms to be group automorphisms. Certain lemmas requiring case analysis will be dealt with in the appendix.

\begin{lemma}
\label{CommutingCycle}
Let $t_1,t_2 \in T$. Then there exists a unique $4$-cycle containing the path $t_2 \to (e) \to t_1$ iff $t_1t_2 = t_2t_1$, and furthermore the cycle will be $(e) \to t_1 \to t_1t_2 \to t_2 \to (e)$.
\end{lemma}

\begin{proof}
See appendix.
\end{proof}

\begin{lemma}
\label{NonCommutingCycle}
Let $t_1,t_2 \in T$ such that $t_1t_2 \neq t_2t_1$. Then the $6$-cycle corresponding to $t_1t_2t_1t_2t_1t_2$ is sent to another cycle of this form under graph automorphisms when $t_1$ and $t_2$ are transpositions. Otherwise, the same statement holds for the $12$-cycle corresponding to $t_1t_2t_1^{-1}t_2^{-1}t_1t_2t_1^{-1}t_2^{-1}t_1t_2t_1^{-1}t_2^{-1}$.
\end{lemma}

\begin{proof}
The case of transpositions was dealt with by Feng \cite{Fen2}. It is easily verified that the latter construction is a cycle when $t_1$ and $t_2$ are not transpositions (it is the union of two cycles when they are transpositions). Also note that no two consecutive edges correspond to commuting generators, and this property is preserved through graph automorphisms by Lemma~\ref{CommutingCycle}. It is natural to try to prove that this is the only $12$-cycle going through $t_1$ and $t_2$ where no two consecutive edges commute. However, this is false, as shown by the following counterexample: Let $a = (1 \ 2 \ 3 \ 4)$, $b = (1 \ 5 \ 6 \ 7)$, $c = (1 \ 8 \ 9 \ 10)$, $d = (1 \ 11 \ 12 \ 13)$. Then $aba^{-1}b^{-1}aba^{-1}b^{-1}aba^{-1}b^{-1} = abcdcb^{-1}a^{-1}bc^{-1}d^{-1}c^{-1}b^{-1} = (e)$. We say that \textit{edge types are preserved} by an automorphism if, whenver $x_1y_1$ and $x_2y_2$ are edges corresponding to the same element of the generating set, then so are $\phi(x_1y_1)$ and $\phi(x_2y_2)$. If we only allow use of the symbols $a,b,a^{-1},b^{-1}$ and assume that edge types are preserved, then this is indeed the only noncommuting $12$-cycle, as demonstrated in the appendix. This leads to a proof of our theorem in a special case, which we will make use of:

\begin{lemma}
\label{TwoGenerators}
Theorem~\ref{CayleyAutomorphism} holds when $|T|=2$, assuming that edge types are preserved.
\end{lemma}

\begin{proof}
The preceding comments show us that commutators of generators map to commutators of generators. Thus $\phi(a)\phi(b) = \phi(ab)$ for all generators $a,b$, so that $\phi(x)\phi(y) = \phi(xy)$ for all $xy$ by induction. The induction itself is sufficiently non-trivial that we feel obliged to include it, but sufficiently technical that we will relegate it to the appendix, even though it requires no case analysis. We have thus shown that all graph automorphisms fixing $(e)$ are in fact group automorphisms as well. That this implies Theorem~\ref{CayleyAutomorphism} we wait to prove in full generality at the end of this section. $\qed$
\end{proof}

Now for any $a,b \in T$, look at $\Gamma_0 = Cay(S_n,\{a,b\}) \subset \Gamma$. The $12$-cycle described above must lie inside $\Gamma_0$. We wish to show that, for any automorphism $\phi \in Aut(\Gamma)$ fixing $(e)$, $\phi(\Gamma_0) = Cay(S_n,\{\phi(a),\phi(b)\})$, from which it will follow that commutators map to commutators in general, and we will have proved Lemma~\ref{NonCommutingCycle}, whence Theorem~\ref{CayleyAutomorphism} follows from the same arguments as in Lemma~\ref{TwoGenerators}.

We will, in fact, prove a stronger contention, namely that if two edges represent the same group element, then their images also represent the same group element. We first offer an automorphism-invariant criterion for determining whether two adjacent edges represent the same group element of the Cayley graph when $T$ is normal.

\begin{lemma}
\label{SameElement}
Let $x \to y \to z$ be a path in $\Gamma$. Then $xy$ and $yz$ represent the same group element if and only if the number of $4$-cycles going through $xy$ equals the number of $4$-cycles going through $yz$.
\end{lemma}

\begin{proof}
Note that if $xy$ and $yz$ correspond to the same group element, then the number of $4$-cycles going through $xy$ definitely equals the number of $4$-cycles going through $yz$ by Lemma~\ref{CommutingCycle}. (Note that this is true even if $T$ consists of $4$-cycles.) The opposite direction is an easy consequence of the normality condition.
\end{proof}

Now note that, by looking at commutativity of edges, we obtain the incidence structure of $Cyc(T)$. Thus the group elements that each edge corresponds to is uniquely determined by which edge each leaf corresponds to (this is simply a consequence of the fact that a tree is determined by the paths between terminal nodes). Thus, given an edge from $v$ corresponding to a leaf $\lambda$, whose pre-image under $\phi$ is $\lambda_0$, it suffices to prove that any edge from an adjacent vertex $w$ corresponding to $\lambda$ also has pre-image $\lambda_0$. First note that unless $|T|=2$, which has already been dispatched of, all leaves commute. We consider two cases: adjacency between $v$ and $w$ is induced by a leaf, or the adjacency is induced by a non-leaf.

\textbf{Case one:} We may assume that all leaves commute, whence we are done by Lemma~\ref{CommutingCycle}.

\textbf{Case two:} By the normality condition, the group element associated with $vw$ must commute with all but one edge, from which we are again done by Lemma~\ref{CommutingCycle}. Then, noting that leaves are mapped to leaves under any graph automorphism, the final leaf only has one place to go (actually, one could make the argument that there are two places to go -- to itself or to its inverse, but both of these edges correspond to the same group element, which is all that we care about).

This completes our contention, so that we are finally done with Lemma~\ref{NonCommutingCycle}.
\end{proof}

By Lemma~\ref{CommutingCycle}, commutativity of edges is preserved through graph automorphisms. Furthermore, cycles are preserved through graph automorphisms. Thus in particular, $\{\phi(t_1),\phi(t_1)\phi(t_2),\phi(t_2),(e)\}$ must be the image of $\{t_1,t_1t_2,t_2,(e)\}$ if $\phi$ is a graph automorphism fixing $(e)$ and $t_1,t_2 \in T$ commute. This implies that $\phi(t_1)\phi(t_2) = \phi(t_1t_2)$. By the same argument, and using Lemma~\ref{NonCommutingCycle}, $\phi(t_1)\phi(t_2) = \phi(t_1t_2)$ if $t_1,t_2 \in T$ don't commute. Thus $\phi(t_1)\phi(t_2) = \phi(t_1t_2)$ for all $t_1,t_2 \in T$. This implies that $\phi$ is not only a graph automorphism but a group automorphism, by the same argument as in Lemma~\ref{TwoGenerators}.

It follows by abuse of notation that $Aut(Cay(S_n,T))_{(e)} \subset Aut(S_n,T)$, where $G_x$ denotes the stabilizer of $x$ under the action of $G$. But it is well-known that $Aut(S_n) \cong S_n$ (the isomorphism being with the inner automorphism group) for $n \neq 6$ \cite{Fen2}, so that $Aut(S_n,T) \subset Aut(Cay(S_n,T))_{(e)}$ when $n \neq 6$ (it is easily verified that any inner automorphism of $S_n$ preserving $T$ must also preserve incidence in $\Gamma$ and is thus a graph automorphism). Note that $n=6$ only when $k=2$, which has already been dispatched, so the theorem holds for all $n$ that we care about. Since $Aut(Cay(S_n,T)) = R(S_n)Aut(Cay(S_n,T))_{(e)}$ and the two subgroups have trivial intersection, we will have a complete characterization of $Aut(Cay(S_n,T))$ if we can show that $R(S_n)$ is normal. This follows since $R(S_n)$ is closed under conjugation by elements of $Aut(S_n,T)$. Thus we have that $Aut(Cay(S_n,T)) \cong R(S_n) \rtimes Aut(S_n,T)$, as stated.
\end{proof}

\begin{comment}
Though it is always regrettable when a result cannot be proven in full generality, we claim that the normality condition is relatively weak. Indeed, given any set $T$, we can define a \textit{normalization of $T$} to be a new $\mathcal{C}$-graph obtained from $T$ by adding another cycle incident on each leaf of $T$. It is easily verified that this results in a normal set.
\end{comment}

\begin{comment}
Though non-normal generating sets for $S_n$ are too big to test, the characterization works for $A_7$ with the generating set  $T = \{(1 \ 2 \ 3),(1 \ 3 \ 2),(1 \ 4 \ 5),(1 \ 5 \ 4),(1 \ 6 \ 7),(1 \ 7 \ 6)\}$, as the automorphism group has size $120960$ (computed by nauty). In this case the group remains the semi-direct product of $A_7$ and the automorphisms of $S_7$ fixing $T$.
\end{comment}

\section{Quasi-hamiltonicity}

For the sake of future work on the conjecture of Rappaport-Strasser and on Hamiltonicity of directed graphs in general, we generalize the work of Gutin and Yeo to on quasi-hamiltonicity to undirected graphs (see \cite{Gut} for the original paper). We will assume that $R \subset V(\Gamma)$.

\begin{definition}
A \textit{cycle factor} in an undirected graph $\Gamma$ is a subgraph of $\Gamma$ such that every vertex has degree $2$.
\end{definition}

\begin{definition}
Let $QH_1(\Gamma,R) := \{e \in E(\Gamma) \ | \ e \cup R$ is in a cycle factor $\}$. For $k>1$, let $QH_k(\Gamma,R) := \{e \in E(\Gamma) \ | \ QH_{k-1}(\Gamma,e \cup R)$ is connected $\}$. Then we say that $\Gamma$ is \textit{$k$-quasi-hamiltonian} if $QH_k(\Gamma,\{\})$ is connected.
\end{definition}

Obviously $k$-quasi-hamiltonicity in an undirected graph implies $k$-quasi-hamiltonicity in the associated digraph. Indeed, it is equivalent to $k$-quasi-hamiltoncity for digraphs if we disallow cycles of length $2$ in the cycle factor. In particular, an undirected graph is Hamiltonian iff it is $(n-2)$-quasi-hamiltonian, since this implies the existence of a cycle factor containing $n-2$ connected vertices (so the last two vertices must also be connected).

\begin{theorem}
\label{qhflow}
Given $\Gamma$, define the bipartite graph $B(\Gamma)$ to have vertex set $T_1 = \{x_1,\ldots,x_m\} \cup T_2 = \{y_1,\ldots,y_m\}$, where $m = |V(\Gamma)|$, and there exists a directed edge from $x_i$ to $x_j$ iff vertices $i$ and $j$ are adjacent in $\Gamma$. Create a flow network where each edge in $B(\Gamma)$ has capacity $1$ and there is a source $s$ with an edge of capacity $2$ into every vertex in $T_1$, similarly an edge of capacity $2$ from every vertex in $T_2$ into a sink $t$. Then there exists a cycle factor in $\Gamma$ containing $R$ iff there exists a flow of $2m$ from $s$ to $t$, such that all edges pertaining to elements of $R$ have flow going through them.
\end{theorem}

\begin{proof}
Suppose that there exists a cycle factor of $\Gamma$ containing $R$. Push flow through $x_iy_j$ and $x_jy_i$ iff the edge $ij$ is in the cycle factor. This gives the desired flow. Now suppose that we have such a flow and wish to construct a cycle factor. It is well-known that we can ``force'' flow to go through an edge by finding an augmenting path containing that edge and then not adding the back-flow through that edge when we push flow through the augmenting path. Thus asking for the existence of such a flow is equivalent to forcing flow through all of the edges pertaining to $R$ (for brevity, from now on we will call this ``forcing flow through $R$'') and asking for the existence of a flow of $2(m-|R|)$ in the resulting graph. Since any choice of augmenting paths must give us the same maximum flow, we can choose any set of augmenting paths that forces flow through $R$. In particular, given any augmenting path $\mathcal{P}$, we can define another path $r(\mathcal{P})$ to be the path obtained by replacing all instances of $x_i$ with $y_i$ (and vice versa) and reversing the orientation of each edge in $\mathcal{P}$. Note that $\mathcal{P}$ and $r(\mathcal{P})$ are edge-disjoint since $\Gamma$ contains no self-loops. If whenever we augment by a path $\mathcal{P}$, we also augment by $r(\mathcal{P})$, then it will be true that $\mathcal{P}$ is an augmenting path iff $r(\mathcal{P})$ is an augmenting path. In particular, we do this while forcing flow through $R$. We then continue to do this while performing the maxflow algorithm. By the symmetry of our algorithm, after we have completed it there will be flow through an edge $x_iy_j$ iff there is flow through an edge $x_jy_i$. Now take the subgraph $\Gamma_0$ of $\Gamma$ formed by all edges $ij$ such that there is flow through $x_iy_j$ in $B(\Gamma)$. Since we have a flow of $2m$ by assumption, every vertex in $\Gamma_0$ has degree $2$, thus is a cycle factor, completing the theorem. $\qed$
\end{proof}

Due to the high degree of symmetry of Cayley graphs, if $\Gamma$ does not contain a Hamiltonian cycle then it is (intuitively) likely to have a quasi-hamiltonicity number sufficiently high that it is infeasible to check. We would thus like a more efficient block to Hamiltonicity for Cayley graphs.

\begin{definition}
A subset $T$ of a group $G$ is said to have a \textit{left coset partition} if there exists a set $S$ such that $s_1T$ and $s_2T$ are disjoint for distinct $s_1,s_2 \in S$, and such that $ST = G$.
\end{definition}

\begin{definition}
A cycle factor is said to be \textit{symmetric} if it is also a left coset partition.
\end{definition}

Note that any Hamiltonian cycle is also a symmetric cycle factor. We can thus define the analogous form of quasi-hamiltonicity where all cycle factors are required to be symmetric. Given a sufficiently crisp characterization of sets with coset partitions, it seems likely that a more effective algorithm for Hamiltonicity blocks could be designed.

\section{Conclusion and Open Problems}

A minor but interesting detail of this paper is the dependence of our bound on $X_0(A)$ on the existence of certain primes. There is no reason to believe that this bound should be strict, and so a more complete understanding of $\mathcal{C}$-trees may be reached by a more precise study of the properties of $X_0(A)$. If the bound is given by explicit constructions, then the result of such a study would also be smaller Cayley graphs that would be more feasible to analyze empirically.

Disregarding our poor understanding of $X_0(A)$, $\mathcal{C}$-trees have been effectively characterized. With this stepping stone, it would be useful to define some more concepts related to $\mathcal{C}$-graphs (in a structurally interesting way). After a tree, the next simplest definition to make is that of a cycle. For a possible idea, we will borrow ideas from matroid theory. We call a set \emph{independent} if it is a subset of a tree, and \emph{dependent} otherwise. A \emph{simple cycle} is then a subset of $T$ that is dependent, but whose proper subsets are all independent. Of course, any other algebraic properties of graphs that could carry over to $\mathcal{C}$-graphs would also be interesting. The author believes that an alternate definition of $\mathcal{C}$-trees leading to a nice matroid structure on the power set of $T$ would make all remaining generalizations transparent. The trees under this structure would also most likely lead to even more structured Cayley graphs.

We have fully characterized the automorphism groups of certain $\mathcal{C}$-trees. We would like a generalization, at the very least, of arbitrary split sets consisting of $k$-cycles. We conjecture that Theorem~\ref{CayleyAutomorphism} holds for all such sets. Similarly, we seek a generalization of Feng's theorem regarding the isomorphism between $Aut(Tra(T))$ and $Aut(S_n,T)$ that gives a relation between the automorphism groups of $Cyc(T)$ and $Aut(S_n,T)$.

We turn to the spectral analysis of the Cayley graphs in question. In particular, we look at the Cayley graphs formed when $|T|=2$. Clearly the largest eigenvalue is $4$ in this case, since the Laplacian is a positive semidefinite matrix. It is interesting to note that the second-largest eigenvalue is $1+\sqrt{3}$ when $k=2$ and $1+\sqrt{7}$ when $k=4$ (the latter result was established empyrically). It is therefore very tempting to conjecture that the second-largest eigenvalue is always $1+\sqrt{2k-1}$, but this is unfortunately nonsense since it can never exceed $4$.

Most importantly, this paper points to a deeper connection between Cayley graphs formed by transpositions and by $k$-cycles. This is structurally apparent in the similarities between the two in terms of commutativity and conjugacy, and indicates that more results should generalize to the case of $k$-cycles. For example, in \cite{Akh} the Coxeter representation of the transpositions is used to gain information about the spectrum of the Cayley graphs. A generalization of the notion of Coxeter representation to account for $k$-cycles would likely allow for the generalization of these results.

Finally, we propose a more general mathematical program to understand the nature of Cayley graphs formed by conjugate generating sets in general, which we believe to be a distinguished variety. So far, all theorems regarding these graphs show that the automorphism group is minimal in a certain sense. We propose the task of finding cases when the group is not minimal, but is close to minimal, and analyzing what happens there, when everything should be more transparent. This should point us towards more general results regarding these graphs.

\section{Appendix}

\subsection{A Bound on Dirichlet's Theorem}

\begin{theorem}
\label{Dirichlet}
For each $n>1$, there is a prime of the form $kn+1$ that divides $\Phi_n(n)$. 
\end{theorem}

\begin{proof}
It can be shown that if $p | \Phi_n(j)$, then either $p | n$ or $p \equiv 1 \pmod{n}$. But $\Phi_n(n) \equiv 1 \pmod{n}$, so we must have the second case. Additionally, $\Phi_n(n) = \prod n-\xi$, for each primitive root of unity $\xi$. But 

$$\displaystyle\prod n-\xi = \sqrt{\prod (n-\xi)(n-\bar{\xi})} = \sqrt{\prod n^2+1-2n\cos(\theta)} > \sqrt{\prod n^2+1-2n} = \prod n-1 \geq 1$$

so we must have some prime dividing $\Phi_n(n)$, and we are done.
\end{proof}

\subsection{Case Analysis for $4$-cycles}

We are asking for $a'$ and $b'$ such that $ab'a'b = (e)$, or equivalently $bab'a' = (e)$. Thus (for the supports of $ab$ and $a'b'$ to be the same) $a'b' \in \{ab,ab^{-1},a^{-1}b,a^{-1}b^{-1},ba\}$. $a'$ cannot be $b^{-1}$ and $b'$ cannot be $a^{-1}$ since this would correspond to a path doubling back on itself. Since the product of two split $k$-cycles is a $2k-1$-cycle, $abab$ is a $2k-1$-cycle, and in particular not the identity. $abab^{-1} = a(bab^{-1})$ is the product of two $k$-cycles with different supports, and so is again not the identity. Similar logic holds for $aba^{-1}b=(aba^{-1})b$. $aba^{-1}b^{-1} = (e)$ implies that $ab=ba$, which is what we want. Finally, $abba = (e)$ implies $aabb=(e)$, which is impossible since $aa$ and $bb$ have different supports.

\subsection{Case Analysis for Commutators}

The full details of this argument can be found at \url{http://www.tjhsst.edu/~jsteinha/Cayley.pdf}. The crux of the argument is repeated use of symmetry, which eventually shows that all interesting cases WLOG start with $abab^{-1}$. We then list out all products of four generators such that no two consecutive generators commute or represent the same element. By looking at what compositions of permutations can send $1$ back to $1$, we reduce essentially to $6$ remaining cases, which are easy to check through simple calculations. This completes our case analysis.

\subsection{Induction Argument for Lemma~\ref{TwoGenerators}}

We have already shown that $\phi(a)\phi(b) = \phi(ab)$ for any automorphism $\phi$ of $\Gamma$ fixing $(e)$. We have the following lemma:

\begin{lemma}
\label{AuxiliaryAutomorphism}
If $\phi$ is a (graph) automorphism of $\Gamma$, then so is $\phi_y = \phi(y^{-1})\phi(yx)$.
\end{lemma}

The proof is a routine verification. Now, we wish to show by induction that

$$\phi(t_1t_2\ldots t_n) = \phi(t_1)\phi(t_2)\ldots \phi(t_n)$$

for all $\phi \in \Gamma$. Now note that

$$\phi(t_1t_2\ldots t_n) = \phi(t_1)\phi_{t_1}(t_2\ldots t_n) = \phi(t_1)\phi_{t_1}(t_2)\ldots \phi_{t_1}(t_n) = \phi(t_1)\ldots\phi(t_n)$$

where the equality between the second and third expressions follows by the inductive step. This completes our induction.

\section{Acknowledgements}

The author would like to thank John Dell, Dave Jensen, Alfonso Gracia-Saz, Jim Lawrence, and Brendan McKay for their help, as well as all of the staff of MathCamp 2007 and Thomas Jefferson High School for Science and Technology.

\singlespacing

\end{document}